\documentclass[12pt]{amsart}
\usepackage{amssymb}
\usepackage{amscd}
\usepackage{amsmath}

\calclayout
\def\comment#1{{\marginpar{\footnotesize#1}}}

\numberwithin{equation}{subsection}

\newtheorem{lemma}[subsection]{Lemma}
\newtheorem{proposition}[subsection]{Proposition}

\theoremstyle{definition}

\theoremstyle{remark}

\newtheorem{remark}[subsection]{Remark}

\newtheorem{example}[subsection]{Example}

\def\OmitComment#1end_of_OmitComment{{}}
\def\comment#1{{\marginpar{\footnotesize#1}}}
\newcommand\Name[1]{{\sc #1}}

\newcommand\C{{\mathbb C}}
\newcommand\N{{\mathbb N}}
\newcommand\Z{{\mathbb Z}}
\newcommand\x{{\vec x}}
\title[Linear Congruences]
{Non-Negative Integer Linear Congruences}
\author{John C. Harris}
\address{30 Marlow Avenue,
 Toronto, Ontario,
 Canada~~M4J~3T9}
\email{harrisj@pathcom.com}

\author{David L. Wehlau}
\address{Department of Mathematics and Computer Science,
         Royal Military College\hfil\break \indent
   PO Box 17000 STN Forces,
Kingston, Ontario, Canada~~K7K~7B4}
\email{wehlau@rmc.ca}
\thanks{Research partially supported by NSERC and ARP}
\thanks{\today}


\begin{document}

\begin{abstract}
  We consider the problem of describing all non-negative integer solutions to  a linear
congruence in many variables.  This question may be reduced to solving the
congruence
 $x_1 + 2x_2 + 3x_3 + \ldots + (n-1)x_{n-1} \equiv 0 \pmod n$ where
  $x_i \in {\mathbb N} = \{0,1,2,\ldots\}$.   We
consider the monoid of solutions of this equation and prove a
conjecture of Elashvili concerning the structure of these
solutions.  This yields a simple algorithm for generating most
(conjecturally all) of the high degree indecomposable solutions of the equation.
\end{abstract}

\maketitle

\section{Introduction}
  We consider the problem of finding all non-negative integer solutions to a
   linear congruence
  $$w_1 x_1 + w_2 x_2 + \ldots + w_r x_r \equiv 0 \pmod {n}\\$$
By a non-negative integer solution, we mean a solution
$A=(a_1,a_2,\ldots, a_r)$ with $a_i \in \N:= \{0,1,2,\ldots\}$ for
all $i=1,2,\ldots,r$.

  As one would expect from such a basic question, this problem has a rich history.
   The earliest published discussion of this problem known to the authors was by \Name{Carl W. Strom}
   in 1931 (\cite{St1}).
   %  Strom also gave a complete solution for $n \leq 10$ (\cite{St2}).
    A number of mathematicians have considered this problem.  Notably \Name{Paul Erd\"os}, \Name{Jacques Dixmier},
    \Name{Jean-Paul Nicolas} (\cite{DEN}), \Name{Victor Kac}, \Name{Richard Stanley} (\cite{K})
    and \Name{Alexander Elashvili} (\cite{E}).

  In particular, \Name{Elashvili} performed a number of computer experiments and made a
number of conjectures concerning the structure of the monoid of
solutions.  Here we prove correct one of Elashvili's conjectures.
This allows us to construct most (conjecturally all) of the
``large'' indecomposable solutions by a very simple algorithm.

  Also of interest are the papers \cite{EJ1}, \cite{EJ2}  by \Name{Elashvili}
  and \Name{Jibladze} and \cite{EJP} by \Name{Elashvili}, \Name{Jibladze} and \Name{Pataraia}
  where the ``Hermite reciprocity'' exhibited by the monoid of solutions is examined.

 \section{Preliminaries}

 \OmitComment
   We consider a system of equations over $\N$:
$$\begin{array}{ccccccccr}
  w_{1,1} x_1 & + & w_{1,2} x_2 & + \ldots + & w_{1,r} x_r & \equiv & 0 &\pmod {n_1}\\
  w_{2,1} x_1 & + & w_{2,2} x_2 & + \ldots + & w_{2,r} x_r & \equiv & 0 &\pmod {n_2}&\qquad(\dag)\\
    \vdots    &   &   \vdots    &          &     \vdots    &        &\vdots &   \vdots\\
  w_{s,1} x_1 & + & w_{s,2} x_2 & + \ldots + & w_{s,r} x_r & \equiv & 0 & \pmod {n_s}
\end{array}$$
where $w_{i,j} \in \N$ for all $i,j$.  We wish to find all solutions
$\x = (x_1,\ldots,x_r) \in \N^r$ of the system $(\dag)$.
  We may find all solutions of $(\dag)$ by iteratively solving the $s$ equations.  Thus we
first find the set $S_1$ of all solutions of the first equation:
$M_1:=\{\x\in\N^r \mid w_{1,1}x_1 + w_{1,2}x_2 + \ldots + w_{1,r}x_r \equiv 0\pmod {n_1}\}$.
\comment{need to say more here}
end_of_OmitComment

We take $\N = \{0,1,2,\ldots\}$ and let $n$ be a positive integer.
Consider the linear congruence
\begin{equation}
\label{one equation}
w_1 x_1 + w_2 x_2 + \ldots + w_r x_r \equiv 0 \pmod n
\end{equation}
where $w_1,w_2,\ldots,w_r \in \Z$ and $x_1,x_2,\dots,x_n$ are unknowns.
We want to describe all solutions $A=(a_1,a_2,\ldots,a_{n-1}) \in {\N}^r$ to this congruence.

  Clearly all that matters here is the residue class of the $w_i$ modulo $n$ and thus we may assume
that $0 \leq w_i < n$ for all $i$.  Also if one of the $w_i$ is
divisible by $n$ then the equation imposes no restriction
whatsoever on $x_i$ and thus we will assume that $1 \leq w_i < n$
for all $i$.

  If $w_1=w_2$ then we may replace the single equation (\ref{one equation}) by the pair of equations
$$w_1 y_1 + w_3 x_3 + \ldots + w_r x_r \equiv 0 \pmod n \qquad {\rm and} \qquad
  x_1 + x_2 = y_1.$$

Thus we may assume that the $w_i$ are distinct and so we have reduced to the case where
$\{w_1,\ldots,w_r\}$ is a subset of $\{1,2,...,n-1\}$.  Now we consider
\begin{equation}
\label{eq}
x_1 + 2 x_2 + 3 x_3+ \ldots + (n-1) x_{n-1} \equiv 0 \pmod n
\end{equation}

The solutions to (\ref{one equation}) are the solutions to (\ref{eq}) with
$x_i = 0$ for all $i \notin \{w_1,\ldots,w_r\}$.
Hence to solve our original problem it suffices to find all solutions to
Equation~(\ref{eq}).

\section{Monoid of Solutions}\label{Monoid}
  We let $M$ denote the set of all solutions to  Equation (\ref{eq}),
$$M := \{\x\in\N^{n-1} \mid x_1 + 2 x_2 + \ldots + (n-1)x_{n-1} \equiv 0\pmod {n}\} \ .$$
Clearly $M$ forms a monoid under componentwise addition, i.e., $M$ is closed under this addition
and contains an additive identity, the {\it trivial solution\/} ${\mathbf 0} = (0,0,\ldots,0)$.

    In order to describe all solutions of (\ref{eq}) explicitly we want to find the set of minimal
generators of the monoid $M$.  We denote this set of generators by $IM$.  We say that a non-trivial
solution
$A \in M$ is {\it decomposable\/} if $A$ can be written as non-trivial sum of two other solutions:
$A = B + C$ where $B, C \neq {\mathbf 0}$.  Otherwise we say that $A$ is {\it indecomposable\/}
(also called {\it non-shortenable\/} in the literature).  Thus $IM$ is the set of indecomposable
solutions.

  We define the {\it degree\/} (also called the {\it height\/} in the literature) of a solution
  $A = (a_1,a_2,\ldots,a_{n-1}) \in M$ by
$\deg(A) = a_1 + a_2 + \ldots + a_{n-1}$ and we denote the set of
solutions of degree $k$ by $M(k) := \{A \in M \mid \deg(A)=k\}$.
Similarly, we let $IM(k)$ denote the set of indecomposable
solutions of degree $k$: $IM(k) = IM \cap M(k)$.

Gordan's Lemma \cite{G} states that there are only finitely many indecomposable
solutions, i.e., that $IM$ is finite.  This is also easy to see directly as follows.
The extremal solutions $E_1:=(n,0,\ldots,0)$, $E_2:=(0,n,0,\ldots,0)$, $\ldots$,
$E_{n-1} := (0,0,\ldots,0,n)$ show
that any indecomposable solution, $(a_1,a_2,\ldots,a_n)$ must satisfy $a_i \leq n$ for all $i$.

  In fact, \Name{Emmy Noether} \cite{N} showed that if $A$ is indecomposable then $\deg(A) \leq n$.
Furthermore $A$ is indecomposable with $\deg(A)=n$ if and only if
$A$ is an extremal solution $E_i$ with $\gcd(i,n)=1$.  For a
simple proof of these results see \cite{S}.

  We define the {\it multiplicity\/} of a solution $A$, denoted $m(A)$ by
$$m(A) := \frac{a_1 + 2a_2 + \ldots + (n-1)a_{n-1}}{n} \ \ .$$

\begin{example}
  Consider $n=4$.  Here
$IM = \{A_1 =(4,0,0), A_2=(0,2,0), A_3=(0,0,4), A_4=(1,0,1),
A_5=(2,1,0), A_6=(0,1,2)\}$. The degrees of these solutions are
$4,2,4,2,3,3$ respectively and the multiplicities are
$1,1,3,1,1,2$ respectively.
\end{example}

\section{The Automorphism Group}\label{Aut}
  Let $G := Aut(\Z/n\Z)$.  The order of $G$ is given by $\phi(n)$ where $\phi$ is the Euler
phi function, also called the totient function.  The elements of $G$ may be represented by
the $\phi(n)$ positive integers less than $n$ and relatively prime to $n$.  Each such integer
$g$ induces a permutation, $\sigma = \sigma_g$, of $\{1,2,\ldots,n-1\}$ given
by $\sigma(i) \equiv gi \pmod n$. Let $A=(a_1,a_2,\ldots,a_{n-1}) \in M$, i.e.,
$a_1 + 2a_2 + \ldots + (n-1)a_{n-1} \equiv 0 \pmod n$.  Multiplying this equation by $g$
gives $(g)a_1 + (2g)a_2 + (3g)a_3+\ldots + (gn-g)a_{n-1} \equiv 0 \pmod n$.
Reducing these new coefficients modulo $n$ and reordering this becomes
$a_{\sigma^{-1}(1)} + 2a_{\sigma^{-1}(2)} + \ldots + (n-1) a_{\sigma^{-1}(n-1)} \equiv 0 \pmod n$.
Thus if $A =(a_1,a_2,\ldots,a_{n-1}) \in M$ then
$g \cdot A := (a_{\sigma^{-1}(1)},a_{\sigma^{-1}(2)} ,\ldots, a_{\sigma^{-1}(n-1)}) \in M$

  Since $g \cdot A$ is a permutation of $A$,
the action of $G$ on $M$ preserves degree, and thus $G$ also acts on each $M(k)$ for
$k \in \N$.  Note however that the action does not preserve multiplicities in general.
Furthermore if $g \in G$ and $A=B+C$ is a decomposable solution, then
$g \cdot A = g \cdot B + g \cdot C$ and therefore $G$ preserves $IM$ and each $IM(k)$.

\begin{example}
  Consider $n=9$.  Here $G$ is represented $\{1,2,4,5,7,8\}$ and the corresponding
six permutations of $\Z/9\Z$ are given by
$\sigma_1 = e$, $\sigma_2 = (1,2,4,8,7,5)(3,6)$, $\sigma_4 = \sigma_2^2 = (1,4,7)(2,8,5)(3)(6)$,
$\sigma_5 = \sigma_2^5 = (1,5,7,8,4,2)(3,6)$, $\sigma_7 = \sigma_2^4 = (1,7,4)(2,5,8)(3)(6)$ and
$\sigma_8 = \sigma_2^3 = (1,8)(2,7)(3,6),(4,5)$.  Thus, for example,
$2\cdot(a_1,a_2,a_3,a_4,a_5,a_6,a_7,a_8)=(a_5,a_1,a_6,a_2,a_7,a_3,a_8,a_4)$ and
$4\cdot(a_1,a_2,a_3,a_4,a_5,a_6,a_7,a_8)=(a_7,a_5,a_3,a_1,a_8,a_6,a_4,a_2)$.
\end{example}

  Note that $G$ always contains the element $n-1$ which is of order 2 and which we also
denote by $-1$.  This element induces the permutation $\sigma_{-1}$ which acts via
$-1\cdot(a_1,a_2,\ldots,a_{n-1})=(a_{n-1},a_{n-2},\ldots,a_3,a_2,a_1)$.

  Let $F(n)$ denote the number of indecomposable solutions to Equation~(\ref{eq}),
$F(n) := \#IM$.  \Name{Victor Kac} \cite{K} showed that the number
of minimal generators for the ring of invariants of $SL(2,\C)$
acting on the space of binary forms of degree $d$ exceeds $F(d-2)$
if $d$ is odd. \Name{Kac} credits \Name{Richard Stanley} for
observing that if $A$ is a solution of multiplicity 1 then $A$ is
indecomposable. This follows from the fact that the multiplicity
function $m$ is a homomorphism of monoids from $M$ to $\N$ and 1
is indecomposable in $\N$.  \Name{Kac} also observed that the
extremal solutions $E_i$ (defined in Section~\ref{Monoid} above)
with $\gcd(i,n)=1$ are also indecomposable.   This gave \Name{Kac}
the lower bound $F(n) \geq p(n) + \phi(n) -1$ where $p(n)$ denotes
the number of partitions of $n$.

%  We could improve this lower bound taking into account the action of the %automorphism group, $G$.

  \Name{Jacques Dixmier}, \Name{Paul Erd\"os} and \Name{Jean-Louis Nicholas} studied the function $F(n)$ and
  significantly improved \Name{Kac}'s lower bound (\cite{DEN}).  They were able to prove that
$$\lim_{n \to \infty}\inf F(n) \cdot \left[ \frac{n^{1/2}}{\log n \cdot \log \log n} p(n) \right]^{-1} > 0\ .$$

  It is tempting to think that the $G$-orbits of the multiplicity 1 solutions would comprise all
elements of $IM$.  This is not true however.  Consider $n=6$.  Then $G$ is a group of order
2, $G=\{1,-1\}$.  The solutions
$A_1 = (1,0,1,2,0)$ and $A_2 = -1 \cdot A_1 = (0,2,1,0,1)$ are both indecomposable and both
have multiplicity 2.

  We define the {\it level\/} of a solution $A$, denoted $\ell(A)$, by
$\ell(A) = \min\{m(g(A)) \mid g \in G\}$.

  Note that
$m(A) + m(-1 \cdot A) = \deg(A)$.  This implies
%$2 \sum_{g \in G} m(g \cdot A) = \deg(A)\phi(n)$.
$2 \sum_{B \in G\cdot A} m(B) = \deg(A)\#(G\cdot A)$,
i.e., that the average multiplicity of the elements in the $G$-orbit of $A$
is half the degree of $A$.

\section{Elashvili's conjectures}
  In \cite{E}, \Name{Elashvili} made a number of interesting and deep conjectures concerning the
structure of the solutions to Equation~(\ref{eq}).  In order to
state some of these conjectures we will denote by $p(t)$ the number of partitions of
the integer $t$.  We also use $\lfloor n/2 \rfloor$ to denote the
greatest integer less than or equal to $n/2$ and $\lceil n/2\rceil
:= n - \lfloor n/2 \rfloor$.

   Conjecture 1:  If $A \in IM(k)$ where %is an indecomposable solution with $\deg(A) \geq \lfloor n/2 \rfloor + 2$
$k \geq \lfloor n/2 \rfloor + 2$ then $\ell(A) = 1$.

  Conjecture 2: If $k \geq \lfloor n/2 \rfloor + 2$ then $IM(k)$
consists of $p(n-k)$ orbits under $G$.

  Conjecture 3: If $k \geq \lfloor n/2 \rfloor + 2$ then $IM(k)$
contains exactly $p(n-k)$ orbits of level 1.

 Here we prove Conjecture~3.  Furthermore we will show that
if $k \geq \lceil n/2 \rceil + 1$ then every orbit of level
1 contains exactly one multiplicity 1 element and has size
$\phi(n)$.  Thus if $k \geq \lceil n/2 \rceil + 1$
then $IM(k)$ contains exactly $\phi(n) p(n-k)$ level 1 solutions.

  This gives a very simple and fast algorithm to generate all the level 1
solutions whose degree, $k$, is at least $\lceil n/2 \rceil + 1$ as follows.
For each partition, $n-k = b_1 + b_2 + \cdots + b_s$, of $n-k$
put $b_{s+1} = \cdots = b_k = 0$  and define $c_i := b_i +1$ for $1\leq i \leq k$.
% n - k < k implies s < k
Then define $A$ via $a_i := \#\{j : c_j=i\}$.
This constructs all multiplicity 1 solutions if $k \geq \lceil n/2 \rceil + 1$.
Now use the action of $G$ to generate the $\phi(n)$ solutions in the orbit of each such 
multiplicity 1 solution.

  If Conjecture~2 is true then this algorithm rapidly produces all elements of
$IM(k)$ for $k \geq  \lfloor n/2 \rfloor + 2$.   This is surprising, since without Conjecture~2,
the computations required to generate the elements of $IM(k)$ become increasingly hard as $k$ increases.

 \section{Proof of Conjecture 3}

  Before proceeding further we want to make a change of variables.  Suppose then that
$A \in M(k)$.  We interpret the solution $A$ as a partition of the integer $m(A)n$ into $k$ parts.
This partition consists of $a_1$ 1's, $a_2$ 2's,$ \ldots,$ and $a_{n-1}$ (n-1)'s.
We write this partition as an {\it unordered\/} sequence (or multi-set) of $k$ numbers:
 $$[y_1,y_2,\ldots,y_k] = [\underbrace{1,1,\ldots,1}_{a_1},
\underbrace{2,2,\ldots,2}_{a_2},\ \ldots \ ,\underbrace{(n-1),(n-1),\ldots,(n-1)}_{a_{n-1}}]$$
The integers $y_1,y_2,\ldots,y_k$ with
$1 \leq y_i \leq n-1$ for $1 \leq i \leq k$ are our new variables for describing $A$.
Given $[y_1,y_2,\ldots,y_k]$ we may easily recover $A$ since $a_i := \#\{j \mid y_j = i\}$.

We have $y_1 + y_2 + \ldots + y_k = m(A)n$.

  Notice that the sequence $y_1-1, y_2-1, \ldots, y_k-1$ is a partition of $m(A)n - k$.
 Furthermore, every partition of $m(A)n-k$ arises from a partition of $m(A)n$ into $k$ parts in
this manner.

  The principal advantage of this new description for elements of $M$ is that it makes the
action of $G$ on $M$ more tractable.  To see this let $g \in G$ be a positive integer less
than $n$ and  relatively prime to $n$.
Then $g \cdot [y_1,y_2,\ldots,y_k] = [gy_1 \pmod n, gy_2 \pmod n, \ldots, gy_k \pmod n]$.

  Now we proceed to give our proof of Elashvili's Conjecture 3.

\begin{proposition}\label{prop:quadratic}
  Let $A \in M(k)$ %where $k \geq \lfloor n/2 \rfloor + 2$.
  and let $1 \leq g \leq n-1$ where $g$ is relatively prime to $n$ represent an element of $G$.
  Write $B = g\cdot A$,
  and  $u = m(A)$ and $v = m(B)$.
If $k \geq gu - v$ then $ug^2 - (k+u+v)g + v(n+1) \geq 0$.
\end{proposition}

\begin{proof}
  Write $A = [y_1,y_2,\ldots,y_k]$ where $y_1 \geq y_2 \geq \ldots \geq y_k$.
For each $i$ with $1\leq i \leq k$ we use the division
algorithm to write $g y_i = q_i n + r_i$ where $q_i \in \N$ and $0 \leq r_i < n$.
Then $B = [r_1, r_2, \ldots, r_k]$.  Note that the $r_i$ may fail to be in decreasing order
and also that no $r_i$ can equal 0.

 Now $gun = g (y_1 + y_2 + \ldots + y_k)
     = (q_1 n + r_1) + (q_2 n + r_2) + \ldots + (q_k n + r_k)
     = (q_1 + q_2 + \ldots +q_k)n + (r_1 + r_2 + \ldots + r_k)$
     where $r_1+r_2+\ldots+r_k=vn$.

 Therefore,  $gu = (q_1+q_2+\ldots+q_k) + v$.

Since $y_1 \geq y_2 \geq \ldots \geq y_k$, we have $q_1 \geq q_2 \geq \ldots \geq q_k$.
Therefore from  $gu - v = \sum_{i=1}^k q_i$ we conclude that
$q_i = 0$ for all $i > gu - v$.
Therefore
\begin{eqnarray*}
\sum_{i=1}^{gu-v} g y_i & = & g \sum_{i=1}^{gu-v} [(y_i-1)+1]\\
    &=& g \sum_{i=1}^{ug-v} (y_i-1) + g(gu-v)\\
    &\leq& g \sum_{i=1}^k (y_i-1)+ g(gu-v)\\
    &=& g (un-k) + g^2u-gv
\end{eqnarray*}

  Also
\begin{eqnarray*}
     \sum_{i=1}^{gu-v} g y_i & = & \sum_{i=1}^{gu-v} (q_i n + r_i)\\
       &=& (gu - v)n + \sum_{i=1}^{gu-v} r_i\\
       &\geq& gun - vn + gu - v
\end{eqnarray*}

  Combining these formulae we obtain the desired quadratic condition
  $ug^2 - (k+u+v)g + v(n+1) \geq 0$.
\end{proof}

   Now we specialize to the case $u=v=1$.  Thus we are considering a pair of solutions $A$ and
 $B = g \cdot A$ both of degree $k$
% \geq \lfloor n/2 \rfloor + 2$
 and both of multiplicity 1.

  \begin{lemma}\label{some ones}
   Let $A \in M(k)$ be a solution of multiplicity 1.
   Write $A = [y_1,y_2,\ldots,y_k]$ where $y_1 \geq y_2 \geq \cdots \geq y_k$.
    If  $k \geq \lfloor n/2 \rfloor + 2$ then $y_{k-2}=y_{k-1}=y_k=1$.
    If  $k \geq \lceil n/2 \rceil + 1$ then $y_{k-1}=y_k=1$.
\end{lemma}

 \begin{proof}  First suppose that $k \geq \lfloor n/2 \rfloor + 2$ and
   assume, by way of contradiction, that $y_{k-2} \geq 2$.  Then
   $n = (y_1 + y_2 + \ldots + y_{k-2}) + y_{k-1} + y_k \geq 2(k-2) + 1 + 1
        \geq 2\lfloor n/2 \rfloor+2 \geq n+1$.

  Similarly if $k \geq \lceil n/2 \rceil + 1$ we
   assume, by way of contradiction, that $y_{k-1} \geq 2$.  Then
   $n = (y_1 + y_2 + \ldots + y_{k-1}) + y_k \geq 2(k-1) + 1
        \geq 2(\lceil n/2 \rceil) + 1 \geq n+1$.
 \end{proof}

\begin{proposition}
  Let $A \in M(k)$  be a solution of multiplicity 1 where $k \geq \lceil n/2 \rceil + 1$.
 Then the $G$-orbit of $A$ contains no other element of multiplicity 1.  Furthermore, $G$ acts faithfully
on the orbit of $A$ and thus this orbit contains exactly $\phi(n)$ elements.
\end{proposition}

\begin{proof}
  Let $B = g\cdot A$ for some $g \in G$ and suppose $B$ also has multiplicity 1.
 Lemma~\ref{some ones}  implies that $B = g\cdot A =[r_1,r_2,\ldots,r_{k-2},g,g]$.  Since
$B$ has multiplicity 1, we have $n=r_1+r_2+\ldots+r_{k-2}+g+g \geq 2g+k-2$ and thus
$g \leq (n-k+2)/2 \leq k/2$.
From this we see that the hypothesis $k \geq gu-v$ is satisfied.  %always satisfied when $u=v=1$.
%Since $k \geq \lceil n/2 \rceil +1$
Therefore by Proposition~\ref{prop:quadratic}, $g$ and $k$ must satisfy the quadratic condition
$$g^2 - (k+2)g + (n+1) \geq 0\ \ .$$

 Let $f$ denote the real valued function $f(g) = g^2-(k+2)g + (n+1)$.
Then $f(1) =n-k \geq 0$ and $f(2) = n+1 -2k < 0$ and thus $f$ has a root in the
interval [1,2).  Since the sum of the two roots of $f$ is $k+2$ we see that the other
root of $f$ lies in the interval $(k,k+1]$.  Thus our quadratic condition implies that
either $g \leq 1$ or else $g \geq k+1$.   But we have already seen that $g \leq k/2$ and thus
we must have $g=1$ and so $A=B$.

  %Thus we have shown that if $A$ is a multiplicity 1 solution and $\deg(A) \geq \lceil n/2 \rceil + 1$ then
This shows that the $G$-orbit of $A$ contains no other element of multiplicity 1.  Furthermore, $G$ acts faithfully
on this orbit and thus it contains exactly $\phi(n)$ elements.
\end{proof}

%\comment{emphasize that level 1 corresponds to a partition}
%\comment{summarize the proof of conjecture 3}

\begin{remark}  Of course the quadratic condition $ug^2 - (k+u+v)g + v(n+1) \geq 0$ can be applied
to cases other than $u=v=1$.  For example, taking $u=v=2$ one can show that a solution of degree $k$
(and level 2) with $k \geq (2n+8)/3$ must have an orbit of size $\phi(n)$ or $\phi(n)/2$.
\end{remark}

%%%%%%%%%%%%%%%%%%%%%%%%%%%%%%%%%%%%%%%%%%%%%%%%%%%%%%%%%%%%
%\newcommand{\noopsort}[1]{} \newcommand{\printfirst}[2]{#1}
%  \newcommand{\singleletter}[1]{#1} \newcommand{\switchargs}[2]{#2#1}
%\providecommand{\bysame}{\leavevmode\hbox to3em{\hrulefill}\thinspace}


\begin{thebibliography}{Ric82b}

\bibitem[DEN]{DEN}Jacques Dixmier, Paul Erd\"os and Jean-Louis Nicolas,
\emph{Sur le nombre d'invariants fondamentaux des formes binaires.
(French) [On the number of fundamental invariants of binary
forms]}, C. R. Acad. Sci. Paris S?r. I Math. \textbf{305} (1987),
no. 8, 319--322.

\bibitem[E]{E}A. Elashvili, \emph{Private Communication}, 1994.

\bibitem[EJ1]{EJ1}A. Elashvili and M. Jibladze,
\emph{Hermite reciprocity for the regular representations of
cyclic groups}, Indag. Math. (N.S.) \textbf{9} (1998), no. 2,
233--238.

\bibitem[EJ2]{EJ2}A. Elashvili and M. Jibladze,
\emph{``Hermite reciprocity'' for semi-invariants in the regular
representations of cyclic groups}, Proc. A. Razmadze Math. Inst.
\textbf{119} (1999), 21--24.

\bibitem[EJP]{EJP}A. Elashvili, M. Jibladze and D. Pataraia,
\emph{Combinatorics of necklaces and ``Hermite reciprocity''},
J. Algebraic Combin. \textbf{10} (1999), no. 2, 173--188.

\bibitem[G]{G} P. Gordan,
\emph{\"Uber die Aufl\"osung linearer Gleichungen mit reellen Coefficienten}, Math. Ann. \textbf{6} (1873), 23--28.

\bibitem[K]{K}Victor G. Kac, \emph{Root systems, representations
 of quivers and invariant theory}, Invariant theory
(Montecatini, 1982), 74--108, Lecture Notes in Math.,
\textbf{996}, Springer, Berlin, 1983.

\bibitem[N]{N}E.~Noether \emph{Der endlichkeitssatz der Invarianten
endlicher Gruppen}, Math. Ann. \textbf{77} (1916) 89--92.

\bibitem[S]{S} B. Schmid, \emph{Finite Groups and Invariant Theory},
Topics in Invariant Theory ( M.-P.~Malliavin Editor), 35--66,
Lecture Notes in Math., \textbf{1478}, Springer-Verlag, Berlin Heidelberg New York, 1991.

\bibitem[St1]{St1} Carl W. Strom, \emph{On complete systems under certain finite
groups}, Bull. Amer. Math. Soc. {\bf 37} (1931) 570--574.

\bibitem[St2]{St2} Carl W. Strom,
\emph{Complete systems of invariants of the cyclic groups of equal
order and degree}, Proc. Iowa Acad. Sci. {\bf 55}, (1948)
287--290.
\end{thebibliography}
\end{document}